\documentclass{gtart}

\def\ifplaintex{\expandafter\ifx\csname documentclass\endcsname\relax}

\def\gtp{{\mathsurround=0pt\it $\cal G\mskip-2mu$eometry \&\ 
$\cal T\!\!$opology $\cal P\!$ublications}}  

\def\recd{{\small Received:\qua\receiveddate\ifx\reviseddate\relax
\else\qquad Revised:\qua\reviseddate\fi\par}} 

\def\lognumber#1{\def\thelognumber{#1}}
\def\volumenumber#1{\def\thevolumenumber{#1}}
\def\volumeyear#1{\def\thevolumeyear{#1}}
\def\papernumber#1{\def\thepapernumber{#1}}
\def\pagenumbers#1#2{\def\startpage{#1}\def\finishpage{#2}}
\def\published#1{\def\publishdate{#1}}

\def\received#1{\def\receiveddate{#1}}
\def\revised#1{\def\reviseddate{#1}}
\def\accepted#1{\def\accepteddate{#1}}

\def\asciiaddress#1{\def\theasciiaddress{#1}}

\long\def\asciiabstract#1{\long\def\theasciiabstract{#1}}

\let\\\par\let\thelognumber\relax\let\thevolumenumber\relax
\let\thepapernumber\relax\let\thevolumeyear\relax\let\startpage\relax
\let\finishpage\relax\let\publishdate\relax\let\receiveddate\relax
\let\reviseddate\relax\let\accepteddate\relax\let\theasciititle\relax
\let\theasciiauthors\relax\let\theasciiaddress\relax
\let\theasciiabstract\relax

\let\theasciiemail\relax

\ifplaintex
\font\logobig=cmssbx10 scaled 3836
\font\logomed=cmssbx10 scaled 2557
\else
\font\logobig=cmssbx10 scaled 4200
\font\logomed=cmssbx10 scaled 2800
\fi

\long\def\makeagttitle{   
\count0=\startpage
\agt\hfill      
\hbox to 45truept{\vbox to 0pt{\vglue -13truept{\logomed A\kern -.37em{\logobig 
T}\kern -.38em G}\vss}\hss}
\break
{\small Volume \thevolumenumber\ (\thevolumeyear)
\startpage--\finishpage\nl
Published: \publishdate}

\vglue .25truein

{\parskip=0pt\leftskip 0pt plus
1fil\def\\{\par\smallskip}{\Large\bf\thetitle}\par\medskip} \vglue
0.05truein

{\parskip=0pt\leftskip 0pt plus 1fil\def\\{\par}{\sc\theauthors}
\par\medskip}%
 
\vglue 0.03truein 

{\small\leftskip 25truept\rightskip 25truept{\bf Abstract}\stdspace\theabstract

{\bf AMS Classification}\stdspace\theprimaryclass
\ifx\thesecondaryclass\relax\else; \thesecondaryclass\fi\par
{\bf Keywords}\stdspace \thekeywords\par}\vglue 7truept

}   

\ifplaintex
\font\phead=cmsl9 scaled 950
\font\pnum=cmbx10 scaled 913
\font\pfoot=cmsl9 scaled 950
\headline{\vbox to 0pt{\vskip -4.5mm\line{\small\phead\ifnum
\count0=\startpage ISSN 1472-2739 (on-line) 1472-2747 (printed)
\hfill {\pnum\folio}\else\ifodd\count0\def\\{ }%
\ifx\theshorttitle\relax\thetitle\else\theshorttitle\fi\hfill{\pnum\folio}
\else\def\\{ and }{\pnum\folio}\hfill\ifx\theshortauthors\relax\theauthors
\else\theshortauthors\fi\fi\fi}\vss}}
\footline{\vbox to 0pt{\vglue 0mm\line{\small\pfoot\ifnum\count0=\startpage
\copyright\ \gtp\hfill\else
\agt, Volume \thevolumenumber\ (\thevolumeyear)\hfill\fi}\vss}}
\else
\font=cmsl9 scaled 1050
\font\lnum=cmbx10 
\font=cmsl9 scaled 1050
\makeatletter
\def\@oddhead{{\small\lhead\ifnum\count0=\startpage ISSN 1472-2739 
(on-line) 1472-2747 (printed)\hfill {\lnum\number\count0}\else\ifodd\count0
\def\\{ }\ifx\theshorttitle\relax \thetitle \else\theshorttitle\fi\hfill
{\lnum\number\count0}\else\def\\{ and }{\lnum\number\count0}
\hfill\ifx\theshortauthors\relax 
\theauthors\else\theshortauthors\fi\fi\fi}}\def\@evenhead{\@oddhead}
\def\@oddfoot{\small\lfoot\ifnum\count0=\startpage\copyright\ \gtp\hfill\else
\agt, Volume \thevolumenumber\ (\thevolumeyear)\hfill\fi}
\def\@evenfoot{\@oddfoot}
\makeatother
\fi
\let\maketitlepage\makeagttitle
\let\makeshorttitle\maketitlepage
\let\maketitle\maketitlepage

\newwrite\gtoutfile
\long\gdef\makeheadfile{  
{\def\\{, }\def\s{ }
\immediate\openout\gtoutfile head.xxx
\immediate\write\gtoutfile{To: math@arxiv.org}
\immediate\write\gtoutfile{Subject: put OR rep NNNNN:ppppp}
\immediate\write\gtoutfile{--text follows this line--}
\immediate\write\gtoutfile{Proxy-for: \ifx\theasciiauthors\relax
\theauthors\else\theasciiauthors\fi\s<\ifx\theasciiemail\relax\theemail\else\theasciiemail\fi>}
\immediate\write\gtoutfile{\noexpand\\}
\immediate\write\gtoutfile{Authors: \ifx\theasciiauthors\relax
\theauthors\else\theasciiauthors\fi}
{\def\\{ }\immediate\write\gtoutfile{Title: \ifx\theasciititle\relax
\thetitle\else\theasciititle\fi}}
\immediate\write\gtoutfile{Subj-class: GT or SG, GR etc}
\immediate\write\gtoutfile{MSC-class: \theprimaryclass\ifx\thesecondaryclass\relax\else, \thesecondaryclass\fi}
\immediate\write\gtoutfile{Journal-ref: Algebr. Geom. Topol. \thevolumenumber\s
(\thevolumeyear) \startpage-\finishpage}
\immediate\write\gtoutfile{Comments: Published by Algebraic and
Geometric Topology at}
\immediate\write\gtoutfile{\s\s\s  http://www.maths.warwick.ac.uk/agt/AGTVol\thevolumenumber/agt-\thevolumenumber-\thepapernumber.abs.html}
\immediate\write\gtoutfile{\noexpand\\}
\immediate\write\gtoutfile{}
\ifx\theasciiabstract\relax
\immediate\write\gtoutfile{\theabstract}\else
\immediate\write\gtoutfile{\theasciiabstract}\fi
\immediate\write\gtoutfile{}
\immediate\write\gtoutfile{\noexpand\\}
\immediate\write\gtoutfile{}
\immediate\closeout\gtoutfile}}  

\def\maketitlepage{\makeagttitle\makeheadfile}
\let\makeshorttitle\maketitlepage
\let\maketitle\maketitlepage

\def\ifplaintex{\expandafter\ifx\csname documentclass\endcsname\relax}

\def\gtp{{\mathsurround=0pt\it $\cal G\mskip-2mu$eometry \&\ 
$\cal T\!\!$opology $\cal P\!$ublications}}  

\def\recd{{\small Received:\qua\receiveddate\ifx\reviseddate\relax
\else\qquad Revised:\qua\reviseddate\fi\par}} 

\def\lognumber#1{\def\thelognumber{#1}}
\def\volumenumber#1{\def\thevolumenumber{#1}}
\def\volumeyear#1{\def\thevolumeyear{#1}}
\def\papernumber#1{\def\thepapernumber{#1}}
\def\pagenumbers#1#2{\def\startpage{#1}\def\finishpage{#2}}
\def\published#1{\def\publishdate{#1}}

\def\received#1{\def\receiveddate{#1}}
\def\revised#1{\def\reviseddate{#1}}
\def\accepted#1{\def\accepteddate{#1}}

\def\asciiaddress#1{\def\theasciiaddress{#1}}

\long\def\asciiabstract#1{\long\def\theasciiabstract{#1}}

\let\\\par\let\thelognumber\relax\let\thevolumenumber\relax
\let\thepapernumber\relax\let\thevolumeyear\relax\let\startpage\relax
\let\finishpage\relax\let\publishdate\relax\let\receiveddate\relax
\let\reviseddate\relax\let\accepteddate\relax\let\theasciititle\relax
\let\theasciiauthors\relax\let\theasciiaddress\relax
\let\theasciiabstract\relax

\let\theasciiemail\relax

\ifplaintex
\font\logobig=cmssbx10 scaled 3836
\font\logomed=cmssbx10 scaled 2557
\else
\font\logobig=cmssbx10 scaled 4200
\font\logomed=cmssbx10 scaled 2800
\fi

\long\def\makeagttitle{   
\count0=\startpage
\agt\hfill      
\hbox to 45truept{\vbox to 0pt{\vglue -13truept{\logomed A\kern -.37em{\logobig 
T}\kern -.38em G}\vss}\hss}
\break
{\small Volume \thevolumenumber\ (\thevolumeyear)
\startpage--\finishpage\nl
Published: \publishdate}

\vglue .25truein

{\parskip=0pt\leftskip 0pt plus
1fil\def\\{\par\smallskip}{\Large\bf\thetitle}\par\medskip} \vglue
0.05truein

{\parskip=0pt\leftskip 0pt plus 1fil\def\\{\par}{\sc\theauthors}
\par\medskip}%
 
\vglue 0.03truein 

{\small\leftskip 25truept\rightskip 25truept{\bf Abstract}\stdspace\theabstract

{\bf AMS Classification}\stdspace\theprimaryclass
\ifx\thesecondaryclass\relax\else; \thesecondaryclass\fi\par
{\bf Keywords}\stdspace \thekeywords\par}\vglue 7truept

}   

\ifplaintex
\font\phead=cmsl9 scaled 950
\font\pnum=cmbx10 scaled 913
\font\pfoot=cmsl9 scaled 950
\headline{\vbox to 0pt{\vskip -4.5mm\line{\small\phead\ifnum
\count0=\startpage ISSN 1472-2739 (on-line) 1472-2747 (printed)
\hfill {\pnum\folio}\else\ifodd\count0\def\\{ }%
\ifx\theshorttitle\relax\thetitle\else\theshorttitle\fi\hfill{\pnum\folio}
\else\def\\{ and }{\pnum\folio}\hfill\ifx\theshortauthors\relax\theauthors
\else\theshortauthors\fi\fi\fi}\vss}}
\footline{\vbox to 0pt{\vglue 0mm\line{\small\pfoot\ifnum\count0=\startpage
\copyright\ \gtp\hfill\else
\agt, Volume \thevolumenumber\ (\thevolumeyear)\hfill\fi}\vss}}
\else
\font=cmsl9 scaled 1050
\font\lnum=cmbx10 
\font=cmsl9 scaled 1050
\makeatletter
\def\@oddhead{{\small\lhead\ifnum\count0=\startpage ISSN 1472-2739 
(on-line) 1472-2747 (printed)\hfill {\lnum\number\count0}\else\ifodd\count0
\def\\{ }\ifx\theshorttitle\relax \thetitle \else\theshorttitle\fi\hfill
{\lnum\number\count0}\else\def\\{ and }{\lnum\number\count0}
\hfill\ifx\theshortauthors\relax 
\theauthors\else\theshortauthors\fi\fi\fi}}\def\@evenhead{\@oddhead}
\def\@oddfoot{\small\lfoot\ifnum\count0=\startpage\copyright\ \gtp\hfill\else
\agt, Volume \thevolumenumber\ (\thevolumeyear)\hfill\fi}
\def\@evenfoot{\@oddfoot}
\makeatother
\fi
\let\maketitlepage\makeagttitle
\let\makeshorttitle\maketitlepage
\let\maketitle\maketitlepage

\newwrite\gtoutfile
\long\gdef\makeheadfile{  
{\def\\{, }\def\s{ }
\immediate\openout\gtoutfile head.xxx
\immediate\write\gtoutfile{To: math@arxiv.org}
\immediate\write\gtoutfile{Subject: put OR rep NNNNN:ppppp}
\immediate\write\gtoutfile{--text follows this line--}
\immediate\write\gtoutfile{Proxy-for: \ifx\theasciiauthors\relax
\theauthors\else\theasciiauthors\fi\s<\ifx\theasciiemail\relax\theemail\else\theasciiemail\fi>}
\immediate\write\gtoutfile{\noexpand\\}
\immediate\write\gtoutfile{Authors: \ifx\theasciiauthors\relax
\theauthors\else\theasciiauthors\fi}
{\def\\{ }\immediate\write\gtoutfile{Title: \ifx\theasciititle\relax
\thetitle\else\theasciititle\fi}}
\immediate\write\gtoutfile{Subj-class: GT or SG, GR etc}
\immediate\write\gtoutfile{MSC-class: \theprimaryclass\ifx\thesecondaryclass\relax\else, \thesecondaryclass\fi}
\immediate\write\gtoutfile{Journal-ref: Algebr. Geom. Topol. \thevolumenumber\s
(\thevolumeyear) \startpage-\finishpage}
\immediate\write\gtoutfile{Comments: Published by Algebraic and
Geometric Topology at}
\immediate\write\gtoutfile{\s\s\s  http://www.maths.warwick.ac.uk/agt/AGTVol\thevolumenumber/agt-\thevolumenumber-\thepapernumber.abs.html}
\immediate\write\gtoutfile{\noexpand\\}
\immediate\write\gtoutfile{}
\ifx\theasciiabstract\relax
\immediate\write\gtoutfile{\theabstract}\else
\immediate\write\gtoutfile{\theasciiabstract}\fi
\immediate\write\gtoutfile{}
\immediate\write\gtoutfile{\noexpand\\}
\immediate\write\gtoutfile{}
\immediate\closeout\gtoutfile}}  

\def\maketitlepage{\makeagttitle\makeheadfile}
\let\makeshorttitle\maketitlepage
\let\maketitle\maketitlepage

\def\ifplaintex{\expandafter\ifx\csname documentclass\endcsname\relax}

\def\gtp{{\mathsurround=0pt\it $\cal G\mskip-2mu$eometry \&\ 
$\cal T\!\!$opology $\cal P\!$ublications}}  

\def\recd{{\small Received:\qua\receiveddate\ifx\reviseddate\relax
\else\qquad Revised:\qua\reviseddate\fi\par}} 

\def\lognumber#1{\def\thelognumber{#1}}
\def\volumenumber#1{\def\thevolumenumber{#1}}
\def\volumeyear#1{\def\thevolumeyear{#1}}
\def\papernumber#1{\def\thepapernumber{#1}}
\def\pagenumbers#1#2{\def\startpage{#1}\def\finishpage{#2}}
\def\published#1{\def\publishdate{#1}}

\def\received#1{\def\receiveddate{#1}}
\def\revised#1{\def\reviseddate{#1}}
\def\accepted#1{\def\accepteddate{#1}}

\def\asciiaddress#1{\def\theasciiaddress{#1}}

\long\def\asciiabstract#1{\long\def\theasciiabstract{#1}}

\let\\\par\let\thelognumber\relax\let\thevolumenumber\relax
\let\thepapernumber\relax\let\thevolumeyear\relax\let\startpage\relax
\let\finishpage\relax\let\publishdate\relax\let\receiveddate\relax
\let\reviseddate\relax\let\accepteddate\relax\let\theasciititle\relax
\let\theasciiauthors\relax\let\theasciiaddress\relax
\let\theasciiabstract\relax

\let\theasciiemail\relax

\ifplaintex
\font\logobig=cmssbx10 scaled 3836
\font\logomed=cmssbx10 scaled 2557
\else
\font\logobig=cmssbx10 scaled 4200
\font\logomed=cmssbx10 scaled 2800
\fi

\long\def\makeagttitle{   
\count0=\startpage
\agt\hfill      
\hbox to 45truept{\vbox to 0pt{\vglue -13truept{\logomed A\kern -.37em{\logobig 
T}\kern -.38em G}\vss}\hss}
\break
{\small Volume \thevolumenumber\ (\thevolumeyear)
\startpage--\finishpage\nl
Published: \publishdate}

\vglue .25truein

{\parskip=0pt\leftskip 0pt plus
1fil\def\\{\par\smallskip}{\Large\bf\thetitle}\par\medskip} \vglue
0.05truein

{\parskip=0pt\leftskip 0pt plus 1fil\def\\{\par}{\sc\theauthors}
\par\medskip}%
 
\vglue 0.03truein 

{\small\leftskip 25truept\rightskip 25truept{\bf Abstract}\stdspace\theabstract

{\bf AMS Classification}\stdspace\theprimaryclass
\ifx\thesecondaryclass\relax\else; \thesecondaryclass\fi\par
{\bf Keywords}\stdspace \thekeywords\par}\vglue 7truept

}   

\ifplaintex
\font\phead=cmsl9 scaled 950
\font\pnum=cmbx10 scaled 913
\font\pfoot=cmsl9 scaled 950
\headline{\vbox to 0pt{\vskip -4.5mm\line{\small\phead\ifnum
\count0=\startpage ISSN 1472-2739 (on-line) 1472-2747 (printed)
\hfill {\pnum\folio}\else\ifodd\count0\def\\{ }%
\ifx\theshorttitle\relax\thetitle\else\theshorttitle\fi\hfill{\pnum\folio}
\else\def\\{ and }{\pnum\folio}\hfill\ifx\theshortauthors\relax\theauthors
\else\theshortauthors\fi\fi\fi}\vss}}
\footline{\vbox to 0pt{\vglue 0mm\line{\small\pfoot\ifnum\count0=\startpage
\copyright\ \gtp\hfill\else
\agt, Volume \thevolumenumber\ (\thevolumeyear)\hfill\fi}\vss}}
\else
\font=cmsl9 scaled 1050
\font\lnum=cmbx10 
\font=cmsl9 scaled 1050
\makeatletter
\def\@oddhead{{\small\lhead\ifnum\count0=\startpage ISSN 1472-2739 
(on-line) 1472-2747 (printed)\hfill {\lnum\number\count0}\else\ifodd\count0
\def\\{ }\ifx\theshorttitle\relax \thetitle \else\theshorttitle\fi\hfill
{\lnum\number\count0}\else\def\\{ and }{\lnum\number\count0}
\hfill\ifx\theshortauthors\relax 
\theauthors\else\theshortauthors\fi\fi\fi}}\def\@evenhead{\@oddhead}
\def\@oddfoot{\small\lfoot\ifnum\count0=\startpage\copyright\ \gtp\hfill\else
\agt, Volume \thevolumenumber\ (\thevolumeyear)\hfill\fi}
\def\@evenfoot{\@oddfoot}
\makeatother
\fi
\let\maketitlepage\makeagttitle
\let\makeshorttitle\maketitlepage
\let\maketitle\maketitlepage

\newwrite\gtoutfile
\long\gdef\makeheadfile{  
{\def\\{, }\def\s{ }
\immediate\openout\gtoutfile head.xxx
\immediate\write\gtoutfile{To: math@arxiv.org}
\immediate\write\gtoutfile{Subject: put OR rep NNNNN:ppppp}
\immediate\write\gtoutfile{--text follows this line--}
\immediate\write\gtoutfile{Proxy-for: \ifx\theasciiauthors\relax
\theauthors\else\theasciiauthors\fi\s<\ifx\theasciiemail\relax\theemail\else\theasciiemail\fi>}
\immediate\write\gtoutfile{\noexpand\\}
\immediate\write\gtoutfile{Authors: \ifx\theasciiauthors\relax
\theauthors\else\theasciiauthors\fi}
{\def\\{ }\immediate\write\gtoutfile{Title: \ifx\theasciititle\relax
\thetitle\else\theasciititle\fi}}
\immediate\write\gtoutfile{Subj-class: GT or SG, GR etc}
\immediate\write\gtoutfile{MSC-class: \theprimaryclass\ifx\thesecondaryclass\relax\else, \thesecondaryclass\fi}
\immediate\write\gtoutfile{Journal-ref: Algebr. Geom. Topol. \thevolumenumber\s
(\thevolumeyear) \startpage-\finishpage}
\immediate\write\gtoutfile{Comments: Published by Algebraic and
Geometric Topology at}
\immediate\write\gtoutfile{\s\s\s  http://www.maths.warwick.ac.uk/agt/AGTVol\thevolumenumber/agt-\thevolumenumber-\thepapernumber.abs.html}
\immediate\write\gtoutfile{\noexpand\\}
\immediate\write\gtoutfile{}
\ifx\theasciiabstract\relax
\immediate\write\gtoutfile{\theabstract}\else
\immediate\write\gtoutfile{\theasciiabstract}\fi
\immediate\write\gtoutfile{}
\immediate\write\gtoutfile{\noexpand\\}
\immediate\write\gtoutfile{}
\immediate\closeout\gtoutfile}}  

\def\maketitlepage{\makeagttitle\makeheadfile}
\let\makeshorttitle\maketitlepage
\let\maketitle\maketitlepage

\def\ifplaintex{\expandafter\ifx\csname documentclass\endcsname\relax}

\def\gtp{{\mathsurround=0pt\it $\cal G\mskip-2mu$eometry \&\ 
$\cal T\!\!$opology $\cal P\!$ublications}}  

\def\recd{{\small Received:\qua\receiveddate\ifx\reviseddate\relax
\else\qquad Revised:\qua\reviseddate\fi\par}} 

\def\lognumber#1{\def\thelognumber{#1}}
\def\volumenumber#1{\def\thevolumenumber{#1}}
\def\volumeyear#1{\def\thevolumeyear{#1}}
\def\papernumber#1{\def\thepapernumber{#1}}
\def\pagenumbers#1#2{\def\startpage{#1}\def\finishpage{#2}}
\def\published#1{\def\publishdate{#1}}

\def\received#1{\def\receiveddate{#1}}
\def\revised#1{\def\reviseddate{#1}}
\def\accepted#1{\def\accepteddate{#1}}

\def\asciiaddress#1{\def\theasciiaddress{#1}}

\long\def\asciiabstract#1{\long\def\theasciiabstract{#1}}

\let\\\par\let\thelognumber\relax\let\thevolumenumber\relax
\let\thepapernumber\relax\let\thevolumeyear\relax\let\startpage\relax
\let\finishpage\relax\let\publishdate\relax\let\receiveddate\relax
\let\reviseddate\relax\let\accepteddate\relax\let\theasciititle\relax
\let\theasciiauthors\relax\let\theasciiaddress\relax
\let\theasciiabstract\relax

\let\theasciiemail\relax

\ifplaintex
\font\logobig=cmssbx10 scaled 3836
\font\logomed=cmssbx10 scaled 2557
\else
\font\logobig=cmssbx10 scaled 4200
\font\logomed=cmssbx10 scaled 2800
\fi

\long\def\makeagttitle{   
\count0=\startpage
\agt\hfill      
\hbox to 45truept{\vbox to 0pt{\vglue -13truept{\logomed A\kern -.37em{\logobig 
T}\kern -.38em G}\vss}\hss}
\break
{\small Volume \thevolumenumber\ (\thevolumeyear)
\startpage--\finishpage\nl
Published: \publishdate}

\vglue .25truein

{\parskip=0pt\leftskip 0pt plus
1fil\def\\{\par\smallskip}{\Large\bf\thetitle}\par\medskip} \vglue
0.05truein

{\parskip=0pt\leftskip 0pt plus 1fil\def\\{\par}{\sc\theauthors}
\par\medskip}%
 
\vglue 0.03truein 

{\small\leftskip 25truept\rightskip 25truept{\bf Abstract}\stdspace\theabstract

{\bf AMS Classification}\stdspace\theprimaryclass
\ifx\thesecondaryclass\relax\else; \thesecondaryclass\fi\par
{\bf Keywords}\stdspace \thekeywords\par}\vglue 7truept

}   

\ifplaintex
\font\phead=cmsl9 scaled 950
\font\pnum=cmbx10 scaled 913
\font\pfoot=cmsl9 scaled 950
\headline{\vbox to 0pt{\vskip -4.5mm\line{\small\phead\ifnum
\count0=\startpage ISSN 1472-2739 (on-line) 1472-2747 (printed)
\hfill {\pnum\folio}\else\ifodd\count0\def\\{ }%
\ifx\theshorttitle\relax\thetitle\else\theshorttitle\fi\hfill{\pnum\folio}
\else\def\\{ and }{\pnum\folio}\hfill\ifx\theshortauthors\relax\theauthors
\else\theshortauthors\fi\fi\fi}\vss}}
\footline{\vbox to 0pt{\vglue 0mm\line{\small\pfoot\ifnum\count0=\startpage
\copyright\ \gtp\hfill\else
\agt, Volume \thevolumenumber\ (\thevolumeyear)\hfill\fi}\vss}}
\else
\font=cmsl9 scaled 1050
\font\lnum=cmbx10 
\font=cmsl9 scaled 1050
\makeatletter
\def\@oddhead{{\small\lhead\ifnum\count0=\startpage ISSN 1472-2739 
(on-line) 1472-2747 (printed)\hfill {\lnum\number\count0}\else\ifodd\count0
\def\\{ }\ifx\theshorttitle\relax \thetitle \else\theshorttitle\fi\hfill
{\lnum\number\count0}\else\def\\{ and }{\lnum\number\count0}
\hfill\ifx\theshortauthors\relax 
\theauthors\else\theshortauthors\fi\fi\fi}}\def\@evenhead{\@oddhead}
\def\@oddfoot{\small\lfoot\ifnum\count0=\startpage\copyright\ \gtp\hfill\else
\agt, Volume \thevolumenumber\ (\thevolumeyear)\hfill\fi}
\def\@evenfoot{\@oddfoot}
\makeatother
\fi
\let\maketitlepage\makeagttitle
\let\makeshorttitle\maketitlepage
\let\maketitle\maketitlepage

\newwrite\gtoutfile
\long\gdef\makeheadfile{  
{\def\\{, }\def\s{ }
\immediate\openout\gtoutfile head.xxx
\immediate\write\gtoutfile{To: math@arxiv.org}
\immediate\write\gtoutfile{Subject: put OR rep NNNNN:ppppp}
\immediate\write\gtoutfile{--text follows this line--}
\immediate\write\gtoutfile{Proxy-for: \ifx\theasciiauthors\relax
\theauthors\else\theasciiauthors\fi\s<\ifx\theasciiemail\relax\theemail\else\theasciiemail\fi>}
\immediate\write\gtoutfile{\noexpand\\}
\immediate\write\gtoutfile{Authors: \ifx\theasciiauthors\relax
\theauthors\else\theasciiauthors\fi}
{\def\\{ }\immediate\write\gtoutfile{Title: \ifx\theasciititle\relax
\thetitle\else\theasciititle\fi}}
\immediate\write\gtoutfile{Subj-class: GT or SG, GR etc}
\immediate\write\gtoutfile{MSC-class: \theprimaryclass\ifx\thesecondaryclass\relax\else, \thesecondaryclass\fi}
\immediate\write\gtoutfile{Journal-ref: Algebr. Geom. Topol. \thevolumenumber\s
(\thevolumeyear) \startpage-\finishpage}
\immediate\write\gtoutfile{Comments: Published by Algebraic and
Geometric Topology at}
\immediate\write\gtoutfile{\s\s\s  http://www.maths.warwick.ac.uk/agt/AGTVol\thevolumenumber/agt-\thevolumenumber-\thepapernumber.abs.html}
\immediate\write\gtoutfile{\noexpand\\}
\immediate\write\gtoutfile{}
\ifx\theasciiabstract\relax
\immediate\write\gtoutfile{\theabstract}\else
\immediate\write\gtoutfile{\theasciiabstract}\fi
\immediate\write\gtoutfile{}
\immediate\write\gtoutfile{\noexpand\\}
\immediate\write\gtoutfile{}
\immediate\closeout\gtoutfile}}  

\def\maketitlepage{\makeagttitle\makeheadfile}
\let\makeshorttitle\maketitlepage
\let\maketitle\maketitlepage

\lognumber{22}
\volumenumber{2}
\volumeyear{2002}
\papernumber{22}
\published{30 May 2002}
\pagenumbers{449}{463}
\received{11 December 2001}
\revised{17 May 2002}
\accepted{28 May 2002}

\usepackage{amssymb,amsmath,graphicx}

\numberwithin{equation}{section}
\newtheorem{thm}{Theorem}[section]
\newtheorem{prp}[thm]{Proposition}

\theoremstyle{definition}
\newtheorem{dfn}[thm]{Definition}
\newtheorem{rmk}[thm]{Remark}
         
\newcommand{\C}{{\mathbb{C}}} 
\newcommand{\R}{{\mathbb{R}}} 
\newcommand{\Z}{{\mathbb{Z}}}

\newcommand{\la}{\langle}
\newcommand{\ra}{\rangle}
\newcommand{\pa}{\partial}
\newcommand{\krn}{\operatorname{ker}}
\newcommand{\braid}{\operatorname{braid}}
\newcommand{\Dcu}{\operatorname{Dcu}}
\newcommand{\Ucu}{\operatorname{Ucu}}
\newcommand{\Lcu}{\operatorname{Lcu}}

\newcommand{\Ecr}{\operatorname{Ecr}}
\newcommand{\Ocr}{\operatorname{Ocr}}

\begin{document}

\title{Framed holonomic knots}                    
\authors{Tobias Ekholm\\ Maxime Wolff}                    

\address{Department of Mathematics, Uppsala University\\P.O. Box 480,
751 06 Uppsala, Sweden}
\secondaddress{D\'epartement de Math\'ematiques et Informatique,
Ecole Normale Sup\'erieure de Lyon\\
46 all\'ee d'Italie, 69364 Lyon C\'edex 07, France} 

\asciiaddress{Department of Mathematics, Uppsala University\\P.O. Box 480,
751 06 Uppsala, Sweden\\and\\Departement de Mathematiques et Informatique,
Ecole Normale Superieure de Lyon\\
46 allee d'Italie, 69364 Lyon Cedex 07, France}

\email{tobias@math.uu.se, mwolff@ens-lyon.fr}                     

\begin{abstract}   
A holonomic knot is a knot in 3-space which arises as the
$2$-jet extension of a smooth function on the circle. 
A holonomic knot associated to a generic function is naturally framed by 
the blackboard framing of the knot diagram associated to the
$1$-jet extension of the function. There are two classical invariants
of framed knot diagrams: the Whitney index (rotation number) $W$ and the
self linking number $S$. 

For a framed holonomic knot we show that $W$ is bounded above by the
negative of the braid index of the knot, and that the sum of 
$W$ and $|S|$ is bounded by the negative of the
Euler characteristic of any Seifert surface of the knot.

The invariant $S$ restricted to framed holonomic knots with
$W=m$, is proved to split into $n$, where $n$ is the largest natural
number with $n\le\frac{|m|}{2}$, integer invariants.  
Using this, the framed holonomic isotopy
classification of framed holonomic knots is shown to be more refined
than the regular isotopy classification of their diagrams.
\end{abstract}

\asciiabstract{A holonomic knot is a knot in 3-space which arises as
the 2-jet extension of a smooth function on the circle.  A holonomic
knot associated to a generic function is naturally framed by the
blackboard framing of the knot diagram associated to the 1-jet
extension of the function. There are two classical invariants of
framed knot diagrams: the Whitney index (rotation number) W and the
self linking number S.  For a framed holonomic knot we show that W is
bounded above by the negative of the braid index of the knot, and that
the sum of W and |S| is bounded by the negative of the Euler
characteristic of any Seifert surface of the knot.  The invariant S
restricted to framed holonomic knots with W=m, is proved to split into
n, where n is the largest natural number with 2n < |m|+1, integer
invariants.  Using this, the framed holonomic isotopy classification
of framed holonomic knots is shown to be more refined than the regular
isotopy classification of their diagrams.}

\maketitle 

\section{Introduction}

\subsection{Holonomic knots and framing}\label{1A}
Let $f\colon S^1\to\R$ be a smooth function. The
{\em holonomic plane curve} and {\em holonomic space curve}  
associated to $f$ are the $1$-jet extension $c$ and the $2$-jet
extension $C$ of $f$, respectively. That is,   
$c(t)=(f(t),f'(t))$, and $C(t)=(f(t),f'(t),f''(t))$
where $t\in S^1$, 
and where $(x_0,x_1)$ and $(x_0,x_1,x_2)$ are linear coordinates on
$\R^2$ and $\R^3$, respectively. 

A {\em holonomic knot} is a holonomic space curve which is an
embedding. A {\em framed holonomic knot} is a  
holonomic knot with associated holonomic plane curve which is
an immersion. If $f$ is a function giving rise to a framed holonomic
knot $C$ then $(f'(t),f''(t))\ne 0$, for all $t\in S^1$, and we
consider the constant vector field $\pa_{x_2}$ as a normal vector
field along $C$. 

In the space of smooth functions on the circle, the functions
with associated holonomic space curve being a (framed) holonomic
knot form an open and dense subset. A {\em (framed) holonomic isotopy}
is a continuous $1$-parameter family of (framed) holonomic knots, or
equivalently a continuous path in the space of (framed) holonomic knots. 

Vassiliev \cite{V} introduced holonomic knots and proved that any knot
class (topological isotopy class of knots) has a holonomic
representative and also that there exists a 
natural isomorphism from finite type invariants of topological knots
to finite type invariants of holonomic knots. 

Birman and Wrinkle \cite{BW} showed that two holonomic knots
which are topologically isotopic are in fact holonomically isotopic. 
From a combinatorial point of view this means that the holonomic isotopy
classification of holonomic knots is identical to the isotopy
classification of their diagrams. 
(A {\em knot diagram} is the image of a generic projection of a knot
to a plane in $\R^3$, decorated with over and under crossing
information at its double points. An {\em isotopy of a knot digram} 
is defined to be a sequence of planar isotopies and Reidemeister moves,
see e.g.\ Kauffman \cite{K}.)

\subsection{Whitney index and self linking number}
Following Kauffman \cite{K}, we say that two knot diagrams which can
be deformed into each other by a sequence 
of planar isotopies and, second- and third Reidemeister moves (i.e.\ the
moves the projections of which are self-tangency- and triple point
instances) are called {\em regularly isotopic}. There
are two simple invariants of regular isotopy: 

Fix an orientation of the ambient $\R^3$. This orientation 
together with a fixed orientation of the projection direction
associated to the knot diagram induce an orientation on the projection
plane. The {\em Whitney index} 
$W$ is the tangential degree of the knot diagram viewed as an oriented
regular plane curve in the projection plane. The {\em self linking
number} $S$ is the linking number of a knot $K$ which projects to the
diagram and a copy of $K$ shifted slightly in the projection direction,
computed using the fixed orientation of the ambient $\R^3$. 

The Whitney formula \cite{W} expresses the Whitney index of a generic
regular plane curve as follows. Let $q$ be a point on $C$ such that
$C$ lies on one side of the tangent line of $C$ at $q$. Let
$\mu(q)=\pm 1$ denote the winding number of $C$ with respect to $q'$,
where $q'$ is the point $q$ shifted slightly into the half plane which
contains the curve. The orientation of $C$ and the point $q$ induces
an ordering of the preimages of a self intersection point $p$ of $C$. Let
$\epsilon_p=\pm 1$ be the sign of the orientation of the plane induced
from the tangent vectors of the ordered branches of $C$ intersecting
at $p$. Then   
$$
W(C)=-\sum_p\epsilon(p) +\mu(q). 
$$
The self linking number is the sum of the crossing signs over all
crossings of the diagram. Hence $W+S$ is an odd integer. 

In our study of framed holonomic knots we use the orientation
$dx_0\wedge dx_1\wedge dx_2>0$ and project along the 
$x_2$-axis oriented by $dx_2>0$ to define $W$ and $S$ for framed holonomic
knots. The ranges of these invariants are easily found: 

\begin{prp}\label{prprange}
Let $C$ be a framed holonomic knot. Then $W(C)<0$, and if $W(C)=-1$ then
$S(C)=0$ and $C$ represents the unknot. Let $m\le -2$ and $n\in\Z$ be
such that $m+n$ is odd. Then there exists a framed holonomic knot with
$W(C)=m$ and $S(C)=n$. 
\end{prp}

\noindent
Proposition \ref{prprange} is proved in Subsection \ref{pfprprange}. 

It is more interesting to consider the ranges of $W$ and $S$
restricted to diagrams representing a fixed knot class $K$. 
It is easy to see that for any integers
$m$ and $n$ such that $m+n$ is odd there exists a diagram $D$ which
represents $K$ with $W(D)=m$ and $S(D)=n$. If the domains of $W$ and $S$
are restricted further   
to framed holonomic knots which represent $K$ the situation changes
drastically.  

\begin{thm}\label{thmfixclass}
Let $C$ be a framed holonomic knot representing the knot class
$K$. Then\vspace{-2mm}
\begin{equation}\label{braid}
W(C)\le -\braid(K),
\end{equation}
where $\braid(K)$ is the braid index of $K$, and
\begin{equation}\label{genus}
W(C)+|S(C)|\le 2g(K)-1,
\end{equation}
where $g(K)$ denotes the genus of $K$. 
\end{thm}

\noindent
Theorem \ref{thmfixclass} is proved in Subsection
\ref{pfthmfixclass}. The proof of \eqref{genus} 
uses the Bennequin inequality \cite{B} from the theory of Legendrian 
knots: if $(x,y,z)$ are coordinates on $\R^3$ and $\R^3$ is oriented by
$dx\wedge dy\wedge dz>0$ then this inequality asserts that for the
$xy$-diagram  $\Gamma_{\C}$ of a knot $\Gamma$ which is everywhere
tangent to the field of hyperplanes $\krn(dz-ydx)$ and which
represents the knot class $K$  
\begin{equation}\label{bineq}
S(\Gamma_{\C})+|W(\Gamma_{\C})|\le 2g(K)-1.
\end{equation}  
It is a curious fact that the roles of $S$ and $W$ in \eqref{genus}
and \eqref{bineq} are reversed.  

\subsection{New invariants of framed holonomic knots} Trace
\cite{T} showed that two knot diagrams $D$ and $D'$ are regularly
isotopic if and only if they represent the same knot class,
$W(D)=W(D')$, and $S(D)=S(D')$. The classification problem for framed
holonomic knots resembles the problem of classifying knot diagrams up
to regular isotopy in the following way.
Regular isotopy is knot diagram isotopy without first
Reidemeister moves (the move which projects to a cusp-instance) and
framed holonomic isotopy is holonomic isotopy without the holonomic
first Reidemeister move, see Figure \ref{holRM}. 

\begin{thm}\label{thmsplit}
On the space of framed holonomic knots with Whitney index equal to
$m$, the invariant $S$ splits. 
More precisely, to each framed holonomic knot $C$ with
$W(C)=m$ there is associated $n$, where $n$ is the largest integer
with $n\le\frac{|m|}{2}$, integers $S_1(C),\dots,S_n(C)$, 
which are invariant under framed holonomic isotopy. Moreover,
\begin{equation}\label{split}
S(C)=\sum_{j=1}^{n} S_j(C).
\end{equation}   
\end{thm}
\noindent
The invariants $S_j$ are defined in Definition \ref{dfnSk} and 
Theorem \ref{thmsplit} is proved in Subsection \ref{pfthmsplit}. 

In Section \ref{Ex} we give examples of framed holonomic knots
representing the same knot class, with the
same $W$ and $S$ but which are
distinguished up to framed holonomic isotopy by the invariants $S_j$.
This shows that the classification of framed holonomic knots up to
framed holonomic isotopy is more refined than the regular isotopy
classification of their diagrams. This result should be compared to
the result of Birman and Wrinkle mentioned in Subsection \ref{1A}. 

\subsection{Holonomic regular homotopy}
A {\em holonomic regular homotopy} is a continuous $1$-parameter
family of regular holonomic plane curves. 
\begin{prp}\label{prpreghom}
Two regular holonomic plane curves are holonomically regularly
homotopic if and only if they have the same Whitney index.
\end{prp}
\noindent
Proposition \ref{prpreghom} is proved in Subsection
\ref{pfprpreghom}. If the 
holonomic requirements in Proposition \ref{prpreghom} are removed one
obtains the classical Whitney-Graustein theorem \cite{W}. The proof we
present is independent of this theorem.

\section{Diagrams of holonomic knots and Reidemeister moves}
For the readers convenience, basic facts on the geometry of
diagrams of holonomic knots are presented. For proofs of these facts,
see \cite{V}, Proposition 1.  

Let $f\colon S^1\to\R$ be a generic function. Then the
$x_0x_1$-diagram $c$ of the framed holonomic knot $C$ associated to
$f$ has the following properties:
\begin{itemize}
\item[{\bf P1}]
$c$ is a regular curve and if $p$ is a point on $c$ in the upper
(lower) half plane and $v$ is 
the unit tangent of $c$ at $p$ then $\la v,\pa_{x_0}\ra>0$ 
($\la v,\pa_{x_0}\ra<0$), where $\la\,,\ra$ denotes the standard inner 
product on $\R^2$.
\item[{\bf P2}]
$c$ meets the $x_0$-axis at right angles at a finite number of
points corresponding to the local extrema of $f$. The curvature of $c$
at such a point $p$ does not vanish and if $p$ corresponds to maximum
(minimum) of $f$ the unit tangent of $c$ at $p$ equals $-\pa_{x_1}$
($\pa_{x_1}$). 
\item[{\bf P3}]
The only self intersection points of $c$ are transverse double points
which lie in the region $\{(x_0,x_1)\colon x_1\ne 0\}$. The crossing
number of a double point in the upper (lower) half plane is negative
(positive) with respect to the orientation $dx_0\wedge dx_1\wedge
dx_2>0$.
\end{itemize}    

\begin{figure}[ht!]
\begin{center}
\includegraphics[angle=0, width=8cm]{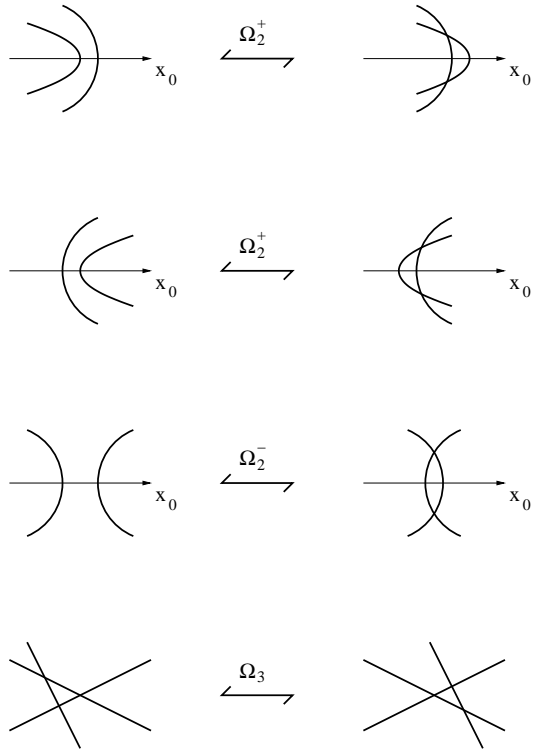}
\end{center}
\caption{Framed holonomic Reidemeister moves}\label{frRM} 
\end{figure}

In generic $1$-parameter families of {\em framed holonomic knots} the
diagram 
changes by planar isotopy which preserve properties {\bf P1}--{\bf 3}
above except for a finite number of instances where one of the
bifurcations in Figure \ref{frRM} occur.
Note that the $\Omega_2$-moves always occur on the $x_0$-axis. The
signs on the $\Omega_2$-moves refer to the signs of the product of 
the second derivatives at the extrema meeting at the self-tangency
moment of the function defining the holonomic knot. The
$\Omega_3$-move depicted occurs either in the upper- or lower half 
plane. 

\begin{figure}[ht!]
\begin{center}\vspace{-2.5cm}
\includegraphics[angle=-90, width=12cm]{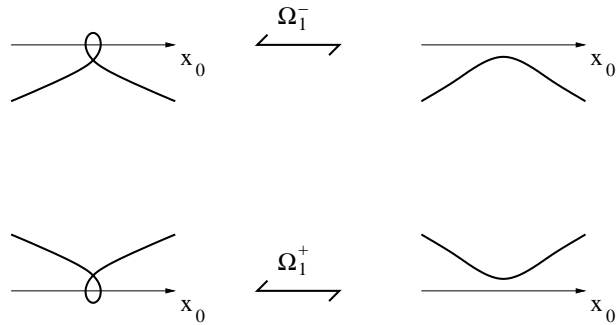}
\vspace{-2.5cm}\end{center}
\caption{Holonomic versions of the first Reidemeister move}\label{holRM} 
\end{figure}

If the word {\em framed} above is omitted the corresponding result is: 
In generic $1$-parameter families of {\em holonomic knots} the diagram
changes by planar isotopy which preserve properties {\bf P1}--{\bf 3}
above except for a finite number of instances where one of the
bifurcations in Figures \ref{frRM} or \ref{holRM} occur.
The signs of the $\Omega_1$-moves in Figure \ref{holRM} refer to the
sign of $x_1$ in the 
half plane where a double point is born or vanishes.

\begin{figure}[ht!]
\begin{center}\vspace{-1.8cm}
\includegraphics[angle=-90, width=12cm]{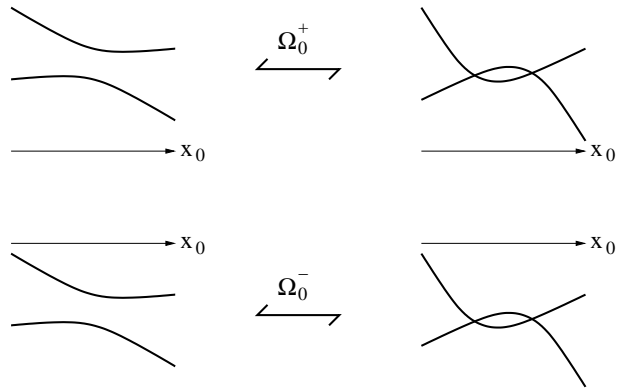}
\vspace{-2.2cm}\end{center}
\caption{Framed holonomic crossing move}\label{CrossM} 
\end{figure}

If we further omit the condition that the holonomic curve be an
embedding then the list of diagram-bifurcations would be further
extended and include also the move in Figure \ref{CrossM} (which might
change the knot class of the holonomic curve).
The signs on the $\Omega_0$-moves in Figure \ref{CrossM} refer to the
half plane where double points are born or vanish.

\subsection{Proof of Proposition \ref{prprange}}\label{pfprprange}
Let $\pi(x_0,x_1,x_2)=(x_0,x_1)$. Let $f$ be a generic function on the 
circle, let $C$ be its associated holonomic knot, and let $c=\pi(C)$.
To see that $W(C)<0$ compute the Whitney index by looking at 
points $p$ on $c$ where the unit tangent equals
$\pa_{x_1}$. These correspond to minima of $f$ and all of them
contribute $-1$ to $W$. The second statement is immediate.

To create a holonomic knot $C$ with $W(C)=m$, $m\le -2$ and $S(C)=n$,
start from the holonomic unknot diagram (which is just the unit
circle). If $n\ge 0$ ($n<0$) apply $\Omega_1^+$ ($\Omega_1^-$) $m-1$
times in such a way that the resulting diagram contains $m-1$ consecutive
loops along the $x_0$-axis. The resulting holonomic knot has $W=m$ and 
$S=-(m-1)$ if $n\ge 0$ ($S=m-1$ if $n<0$). Finally, if $n\ge 0$ apply
$\Omega_0^-$ $\frac{n+m-1}{2}$ times and if $n<0$ apply $\Omega_0^+$ 
$\frac{n+m-1}{2}$ times to create new double points. The resulting
holonomic knot $C$ then has $W(C)=m$ and $S(C)=n$, as desired.\qed 

\subsection{Proof of Proposition \ref{prpreghom}}\label{pfprpreghom}
Let $f$ be a function with associated
holonomic plane curve $c_f$ which is an immersion. If $\phi$ is a
diffeomorphism of $S^1$ then also $f\circ \phi$ gives rise to a regular
holonomic plane curve.

Let $g$ be a  function with regular plane holonomic curve $c_g$
with $W(c_f)=W(c_g)$. Perturb $f$ and $g$ so that they become Morse
functions. Then the proof of Proposition
\ref{prprange} implies that they have the same number of local
extrema. Let $\phi_s$, $0\le s\le 1$ be a diffeotopy of $S^1$ which
moves each critical point of $f$ to a critical point of $g$ of the
same index. Then the critical sets of $\hat g=g\circ \phi_1$ 
and of $f$ agree. Moreover, if $t$ is local maximum (minimum) of $f$
then it is a local maximum 
(minimum) of $g\circ\phi_1$. Let $(s,r)$ be coordinates on the cylinder
$S^1\times\R$ and consider the vector field 
$V(s,r)=(f(s)-\hat g(s))\pa_r$. Let $\Phi^\rho$ be the flow of
$V$. If $\hat g_\rho$ is the function with graph
$\Phi^\rho(\Gamma_{\hat g})$, where $\Gamma_{\hat g}$ is the graph of
$\hat g$, then $\hat g_\rho$ has a regular associated holonomic curve for
each $\rho\le 1$ and $\hat g_1=f$. 

These two deformations together give the desired holonomic regular
homotopy.\qed
         
\section{Holonomic knots and front projections of Legendrian knots}

\subsection{The front and complex projections of a Legendrian knot}
\begin{figure}[ht!]
\begin{center}\vspace{-1.7cm}
\includegraphics[angle=-90, width=8cm]{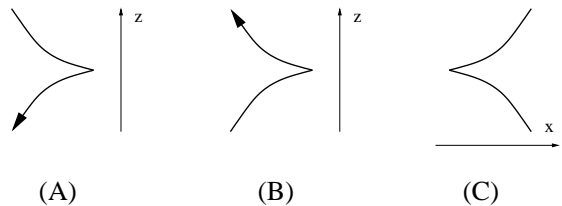}
\vspace{-1.7cm}\end{center}
\caption{(A) Down-cusp, (B) Up-cusp, and (C) left-cusp.}\label{legcusp} 
\end{figure}

Let $\Gamma$ be a knot in $\R^3$ with coordinates $(x,y,z)$
everywhere tangent to the plane field $\{\krn(dz-ydx)\}$. That is,
$\Gamma$ is a {\em Legendrian knot}. Assume moreover that $\Gamma$ is
generic among Legendrian knots, then the projection $\Gamma_F$ of
$\Gamma$ to the $xz$-plane is a curve with transverse double points,
isolated cusps, and without vertical tangencies. Moreover, given any
curve in the $xz$-plane with these properties, there exists a unique
Legendrian knot which projects to this curve. We associate the
following numbers to $\Gamma_F$: 

First we count cusps,
let $\Dcu(\Gamma_F)$, $\Ucu(\Gamma_F)$, and $\Lcu(\Gamma_F)$
denote the number of
down-cusps, up-cusps, and left-cusps  respectively of $\Gamma_F$, 
see Figure \ref{legcusp}. 
 
\begin{figure}[ht!]
\begin{center}\vspace{-2.3cm}
\includegraphics[angle=-90, width=10cm]{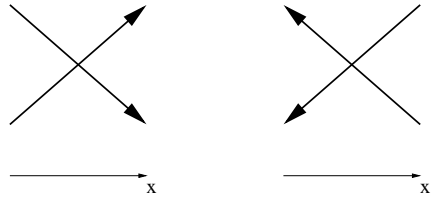}
\vspace{-2.3cm}\end{center}
\caption{Crossing points with tangents with $x$-components of equal
sign}\label{Ecr}  
\end{figure}

\begin{figure}[ht!]
\begin{center}\vspace{-2.3cm}
\includegraphics[angle=-90, width=10cm]{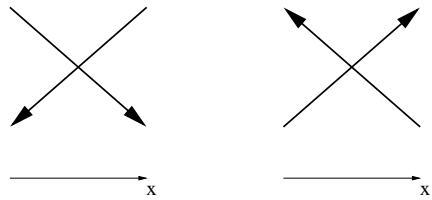}
\vspace{-2.3cm}\end{center}
\caption{Crossing points with tangents with $x$-components of opposite
signs}\label{Ocr}  
\end{figure}

Second we count crossings,
let $\Ecr(\Gamma_F)$ denote the number of crossing points where
the tangent vectors has $x$-components of the same sign and
$\Ocr(\Gamma_F)$ the number of crossing points where
the tangent vectors has $x$-components of the opposite sign. See
Figures \ref{Ecr} and \ref{Ocr}.
 
The projection $\Gamma_{\C}$ of $\Gamma$ to the $xy$-plane is a generic
knot diagram. It is straightforward to check that the invariants
$W(\Gamma_{\C})$, where we use the orientation given by $dx\wedge dy$
in the $xy$-plane, and $S(\Gamma_{\C})$, where we use the orientation
given by $dx\wedge dy\wedge dz$ in space, can be computed 
from data of $\Gamma_F$ as follows, 
\begin{align}
W(\Gamma_{\C})&=\frac12\left(\Dcu(\Gamma_F)-\Ucu(\Gamma_F)\right),\\
S(\Gamma_{\C})&=\Ecr(\Gamma_F)-\Ocr(\Gamma_F)-\Lcu(\Gamma_F).
\end{align}

\subsection{Legendrian knots associated to a holonomic one}\label{legtohol}
Let $C$ be a framed holonomic knot. We associate two
Legendrian knots $\Gamma^+$ and $\Gamma^-$, everywhere
tangent to $\krn(dx_1-x_2dx_0)$, to $C$ 
by describing their front projections (in the $x_0x_1$-plane). The
resulting Legendrian knots lie in $\R^3$ oriented by 
$dx_0\wedge dx_2\wedge dx_1>0$.  

\begin{figure}[ht!]
\begin{center}\vspace{-2cm}
\includegraphics[angle=-90, width=10cm]{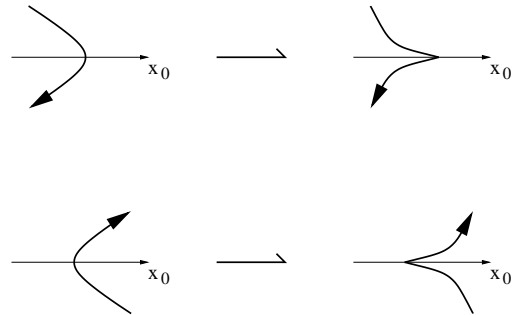}
\vspace{-2cm}\end{center}
\caption{Replacing vertical tangencies with cusps}\label{tang->cusp}  
\end{figure}

\begin{figure}[ht!]
\begin{center}\vspace{-2.7cm}
\includegraphics[angle=-90, width=10cm]{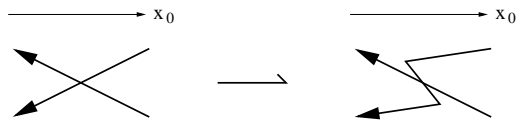}
\vspace{-2.7cm}\end{center}
\caption{Inserting a zig-zag in the lower half plane}\label{zigzagl}  
\end{figure}

The first step in the construction of the fronts of $\Gamma^{+}$ and
$\Gamma^{-}$ is the same in both cases: 

The points where the diagram of
$C$ has vertical tangents are all confined to the $x_0$-axis. 
Replace neighborhoods of such points in
the diagram with cusped arcs as described in Figure \ref{tang->cusp}.
The second step however differs: 

To obtain the front of $\Gamma^+$ we insert a zig-zag as in
Figure \ref{zigzagl} at all crossings in the  lower half plane
and keep the crossings in the upper half plane as they are.

\begin{figure}[ht!]
\begin{center}\vspace{-2.7cm}
\includegraphics[angle=-90, width=10cm]{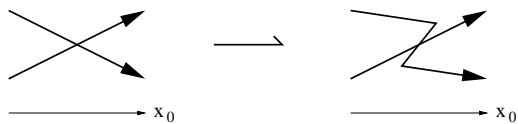}
\vspace{-2.7cm}\end{center}
\caption{Inserting a zig-zag in the upper half plane}\label{zigzagu}  
\end{figure}

To obtain the front of $\Gamma^-$ we insert a zig-zag as in
Figure \ref{zigzagu} at all crossings in the upper half plane and keep
the crossings in the lower half plane as they are.

It is easy to check that $\Gamma^-$ ($\Gamma^+$) is
topologically isotopic to the knot $C$ in $\R^3$ equipped with the
orientation $dx_0\wedge dx_1\wedge dx_2>0$ ($dx_0\wedge dx_1\wedge
dx_2<0$). 

\subsection{Proof of Theorem \ref{thmfixclass}}\label{pfthmfixclass}
Using $\Omega_2$-moves, we may obtain a closed
braid representation of $K$. Equation \eqref{braid} follows.

To prove \eqref{genus}, let $H_+(C)$ ($H_-(C)$) denote the number of
intersection points of $C$ in the upper (lower) half plane. 
Then $S(C)=H_-(C) - H_+(C)$. As noted before, $W(C)$ is the negative of the
number of local minima of the function $f$ giving rise to $C$. 

Let $\Gamma^+$ and $\Gamma^-$ be the Legendrian knots associated to
$C$ as in Subsection \ref{legtohol}. Then
\begin{align*}
&\Dcu(\Gamma^+_F)=-W(C)+2H_-(C),\\
&\Dcu(\Gamma^-_F)=-W(C)+2H_+(C),\\
&\Ucu(\Gamma^+_F)=\Ucu(\Gamma^-_F)=-W(C),
\end{align*}
and hence 
\begin{align}\label{WG+}
W(\Gamma_{\C}^+)&=H_-(C),\\ 
W(\Gamma_{\C}^-)&=H_+(C).\label{WG-}
\end{align}
Also,
\begin{align*}
&\Lcu(\Gamma^+_F)=-W(C)+H_-(C),\\
&\Lcu(\Gamma^-_F)=-W(C)+H_+(C),\\
&\Ecr(\Gamma^+_F)=\Ocr(\Gamma^-_F)=H_+(C),\\
&\Ecr(\Gamma^-_F)=\Ocr(\Gamma^+_F)=H_-(C),
\end{align*}
and hence
\begin{align}\label{SG+}
S(\Gamma_{\C}^+)&=H_+(C) - 2H_-(C) + W(C),\\
S(\Gamma_{\C}^-)&=H_-(C) - 2H_+(C) + W(C).\label{SG-}
\end{align}
Combining \eqref{WG+} and \eqref{SG+}, respectively 
\eqref{WG-} and \eqref{SG-} with the Bennequin inequality
\eqref{bineq} yields  
\begin{align*}
-S(C)+W(C) &\le 2g(K)-1\text{ and }\\
S(C)+W(C) &\le 2g(K)-1,
\end{align*}
since the genus does not depend on the orientation of the ambient
space. The theorem follows.\qed

\section{Splitting the self linking number}
Consider the diagram of a framed holonomic knot $C$. The $x_0$-axis divides
the diagram into cyclically ordered arcs $(X_i,Y_i)$, $i=1,\dots,m$,
where the $X_i$ lies in the upper half plane, the $Y_i$ in the
lower, and where $-m=W(C)$.

Let $(A_i,A_j)=(X_i,X_j)$ or $(A_i,A_j)=(Y_i,Y_j)$ where $i\ne
j$. Define  
\begin{equation*}
\delta(A_i,A_j)=\begin{cases}
                1 &\text{if $\pa A_i$ is contained in an
                unbounded component of $\R-\pa A_j$,}\\
                0 & \text{otherwise.}
                \end{cases}
\end{equation*}
Define 
$$
N(A_i,A_j)=|A_i\cap A_j|+\delta(A_i,A_j),
$$
where $|S|$ denotes the number of points in the set $S$.
 
Let $\tilde X_i$ and $\tilde Y_i$ denote the preimages of $X_i$ and
$Y_i$, for $i=1,\dots,m$. Let $x_i$ and $y_i$ denote the midpoints of
$\tilde X_i$ and $\tilde Y_i$, respectively. 

Consider two arcs $X_i$ and $X_j$, $i\ne j$. Let
$\gamma(x_i,x_j)$ denote the unique oriented arc connecting $x_i$ to
$x_j$ with orientation agreeing with that of the circle. 
Define the {\em cyclic distance of $X_i$ and $X_j$} as 
\begin{equation*}
d(X_i,X_j)=\min\Bigl\{|\gamma(x_i,x_j)\cap\{y_1,\dots,y_m\}|,
|\gamma(x_j,x_i)\cap\{y_1,\dots,y_m\}|\Bigr\}.
\end{equation*} 
Define the cyclic distance of arcs $Y_i$ and $Y_j$ analogously.

\begin{dfn}\label{dfnSk}
Define
$$
S_k(C)=\frac12\left(\sum_{\{(Y_i,Y_j)\colon d(Y_i,Y_j)=k\}}N(Y_i,Y_j)
-\sum_{\{(X_i,X_j)\colon d(X_i,X_j)=k\}}N(X_i,X_j)\right).
$$
\end{dfn}

\begin{rmk}
In terms of defining functions, the terms in the definition of $S_k$
can be interpreted as follows. 
Let $f\colon S^1\to\R$ be a function with associated framed holonomic
knot $C$. Consider $f$ as a periodic function with period $T$ such
that $f(0)=f(T)$ is the global minimum of $f$. Let
$\Gamma_f\subset[0,T]\times\R\subset\R^2$ denote the graph of $f$. 
Then the arcs $X_i$ ($Y_i$) are the holonomic curves
corresponding to restrictions of $f$ to subintervals of $[0,T]$, where
$f$ is increasing (decreasing). If $(x,y)$ are coordinates on $\R^2$
then $|A_i\cap A_j|$ equals the number of lines
$l_a=\{y=a\}$, $a\in\R$ which intersect $A_i$ and $A_j$ at equal
angles, and $\delta(A_i,A_j)=1$ if no $l_a$ intersect both $A_i$ and
$A_j$, otherwise it is $0$.   
\end{rmk}

\subsection{Proof of Theorem \ref{thmsplit}}\label{pfthmsplit}
We  check that $S_k$ is invariant under framed holonomic
Reidemeister moves. For $\Omega_3$ this is immediate.

An $\Omega_2^+$-move involving two arcs $X$ and $X'$ with
distance $d(X,X')=k$ involves also two arcs $Y$ and $Y'$ with
$d(Y,Y')=k$. At the move $\delta(X,X')$ and $\delta(Y,Y')$ are
unchanged and the change in $|X\cap X'|$ and $|Y\cap Y'|$ are the
same. Hence $S_k$ remains constant.

At an $\Omega_2^-$-move involving arcs $X$ and $X'$, $Y$ and $Y'$ 
the change in $|X\cap X'|$ and $|Y\cap Y'|$ equals the change in
$\delta(X,X')$ and $\delta(Y,Y')$, respectively. Hence $S_k$ remains
constant.  

To prove \eqref{split} note that by using the $\Omega_2$-moves we may
move any framed holonomic knot diagram in such a way that its diagram
is a closed braid with braid-axis parallel to the
$x_2$-direction. (The linking number of this axis oriented in the
positive $x_2$-direction and the holonomic knot with its natural
orientation is negative.) Under such deformations both $S$ and
$S_1,\dots,S_n$, remain constant. Moreover for a diagram which is a
closed braid $\delta(X_i,X_j)=\delta(Y_i,Y_j)=0$ for all
$i,j$. Hence both the left and right hand sides of \eqref{split} are
equal to the difference of the number of double points in the lower
and upper half planes. \qed   

\section{Examples}\label{Ex}

\begin{figure}[ht!]
\begin{center}\vspace{-1cm}
\includegraphics[angle=-90, width=6cm]{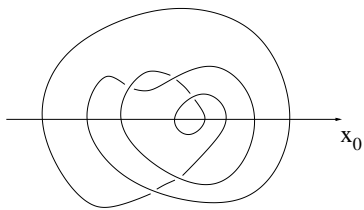}
\vspace{-0.8cm}\end{center}
\caption{The framed holonomic knot $K_1$}\label{knot1}  
\end{figure}

\begin{figure}[ht!]
\begin{center}\vspace{-0.8cm}
\includegraphics[angle=-90, width=6cm]{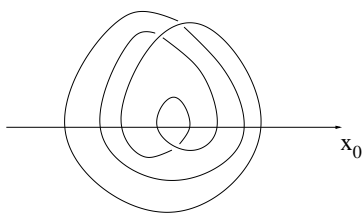}
\vspace{-0.8cm}\end{center}
\caption{The framed holonomic knot $K_2$}\label{knot2}  
\end{figure}

The framed holonomic knots $K_1$ in Figure \ref{knot1} and $K_2$ in
Figure \ref{knot2} both represent the unknot. Since $S(K_1)=S(K_2)=-1$
and $W(K_1)=W(K_2)=-4$, $K_1$ and $K_2$ are regularly isotopic. 
Since $K_1$ is a closed braid $\delta(X_i,X_j)=0=\delta(Y_i,Y_j)$ for
all $i,j$. Noting that all three intersection points in the diagram of
$K_1$ in the upper half plane are intersections between arcs of cyclic
distance $1$, and that the two intersection points in the lower half plane
are intersections of arcs of cyclic distance $1$ respectively $2$, we
conclude that $S_1(K_1)=-2$ and $S_2(K_1)=1$. A similar calculation
gives $S_1(K_2)=0$ and $S_2(K_2)=-1$. Hence $K_1$ and $K_2$ are {\em
not} framed holonomically isotopic. 

\begin{figure}[ht!]
\begin{center}\vspace{-1cm}
\includegraphics[angle=-90, width=6cm]{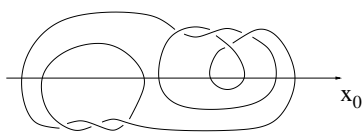}
\vspace{-1cm}\end{center}
\caption{The framed holonomic knot $K_3$}\label{knot3}  
\end{figure}

\begin{figure}[ht!]
\begin{center}\vspace{-1cm}
\includegraphics[angle=-90, width=6cm]{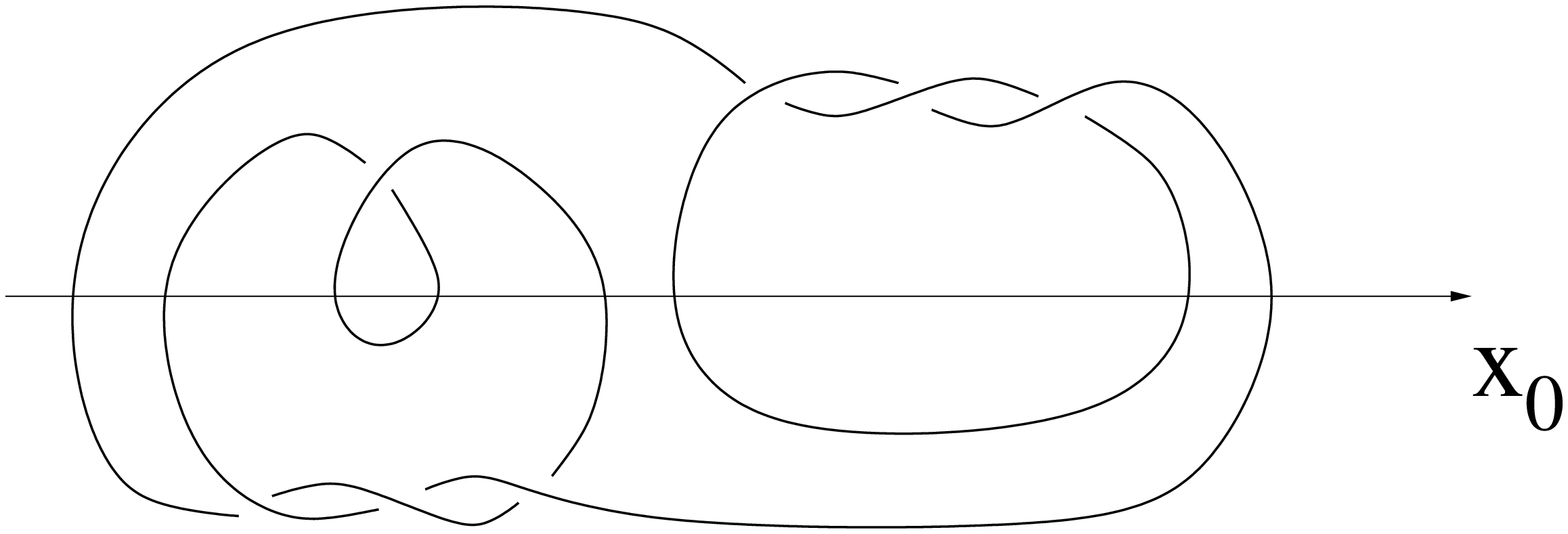}
\vspace{-1cm}\end{center}
\caption{The framed holonomic knot $K_4$}\label{knot4}  
\end{figure}

The framed holonomic knots $K_3$ in Figure \ref{knot3} and $K_4$ in
Figure \ref{knot4} both represent the connected sum of the trefoil
and its mirror image. Since $S(K_3)=S(K_4)=-1$
and $W(K_3)=W(K_4)=-4$, $K_3$ and $K_4$ are regularly isotopic. However,
$S_1(K_3)=0$ and $S_2(K_3)=-1$ but
$S_1(K_4)=-4$ and $S_2(K_4)=3$ so $K_3$ and $K_4$ are {\em not}
framed holonomically isotopic.

%
%
%
%

\Addresses\recd

\end{document}